\documentclass[preprint,12pt,3p,sort]{elsarticle}
\usepackage{amssymb}
\usepackage[colorlinks,bookmarksopen,bookmarksnumbered,citecolor=red,urlcolor=red]{hyperref}
\usepackage{moreverb,dblfloatfix}
\usepackage{amsmath}
\usepackage{algorithmic}
\usepackage{algorithm}
\usepackage{amssymb}
\usepackage{graphicx,color}
\usepackage{float}	
\usepackage{amsmath,amssymb,soul}
\usepackage{subfigure}
\usepackage{rotating}
\usepackage{multirow}
\usepackage{placeins}
\usepackage[normalem]{ulem}
\usepackage{appendix}
\usepackage{listings}
\usepackage{amstext}
\usepackage{epsfig}
\usepackage{fancyhdr}
\usepackage{url}
\usepackage{amsfonts,bm}
\usepackage{upgreek}
\usepackage{setspace}
\usepackage{float}
\usepackage{tabularx}
\usepackage{xcolor}
\usepackage{graphics}
\usepackage{subfigure,epsfig}
\usepackage{mhchem,ulem}
\usepackage{epstopdf}
\usepackage[numbers]{natbib}
\usepackage{soul}
\usepackage{dsfont}

\newcommand{\bdelta}{ {\boldsymbol\delta} }

\renewcommand{\Re}{\mathds{R}}

\makeatletter
\def\ps@pprintTitle{%
 \let\@oddhead\@empty
 \let\@evenhead\@empty
 \def\@oddfoot{}%
 \let\@evenfoot\@oddfoot}
\makeatother

\begin{document}

\begin{frontmatter}

\title{A Learning Based Approach for Uncertainty Analysis \\ in Numerical Weather Prediction Models}

 \author[labela]{Azam Moosavi} 
 \author[labelb]{Vishwas Rao} 
 \author[labela]{Adrian Sandu} 
 \address[labela]{Computational Science Laboratory, Department of Computer Science   \\
Virginia Polytechnic Institute and State University, Blacksburg, VA 24060, USA \\
E-mail: \{azmosavi,\,asandu7\,\}@vt.edu}
 \address[labelb]{Mathematics and Computer Science Division, Argonne National Laboratory, Lemont, IL, USA \\
E-mail: vhebbur@anl.gov}
\thispagestyle{empty}
\setcounter{page}{0}

\begin{Huge}
\begin{center}
Computer Science Technical Report CSTR-{\tt2018-2} \\
\end{center}
\end{Huge}
\vfil
\begin{huge}
\begin{center}
Azam Moosavi, Vishwas Rao, Adrian Sandu
\end{center}
\end{huge}

\vfil
\begin{huge}
\begin{it}
\begin{center}
``{\tt A Learning Based Approach for Uncertainty Analysis in Numerical Weather Prediction Models}''
\end{center}
\end{it}
\end{huge}
\vfil

\begin{large}
\begin{center}
Computational Science Laboratory \\
Computer Science Department \\
Virginia Polytechnic Institute and State University \\
Blacksburg, VA 24060 \\
Phone: (540)-231-2193 \\
Fax: (540)-231-6075 \\ 
Email: \url{azmosavi@cs.vt.edu} \\
Web: \url{http://csl.cs.vt.edu}
\end{center}
\end{large}

\vspace*{1cm}

\begin{tabular}{ccc}
\includegraphics[width=2.5in]{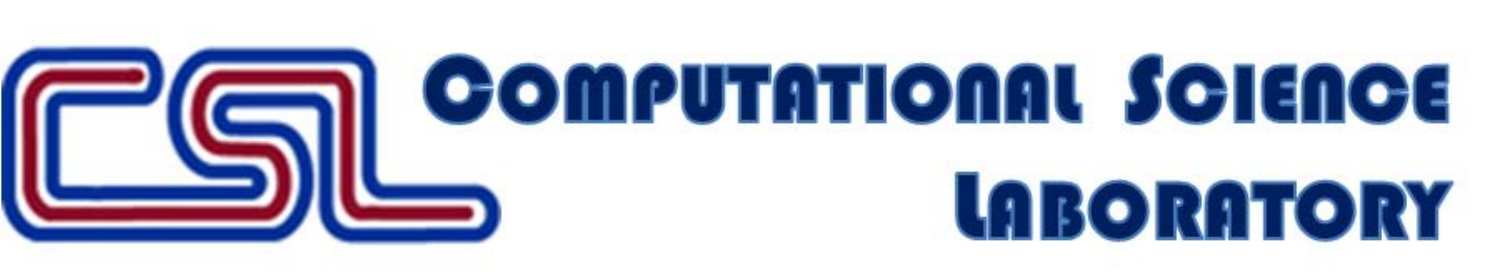}
&\hspace{2.5in}&
\includegraphics[width=2.5in]{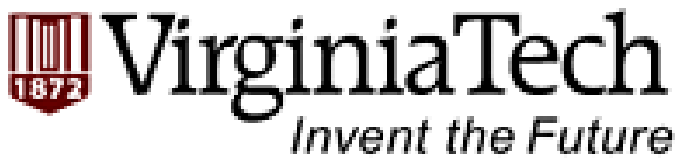} \\
{\bf\em Innovative Computational Solutions} &&\\
\end{tabular}

\newpage

\begin{abstract}
Complex numerical weather prediction models incorporate a variety of physical processes, each described by multiple alternative physical schemes with specific parameters.  The selection of the physical schemes and the choice of the corresponding physical parameters during model configuration can significantly impact the accuracy of model forecasts. There is no combination of physical schemes that works best for all times, at all locations, and under all conditions.  It is therefore of considerable interest to understand the interplay between the choice of physics and the accuracy of the resulting forecasts under different conditions. 

This paper demonstrates the use of machine learning techniques to study the uncertainty in numerical weather prediction models due to the interaction of multiple physical processes.
The first problem addressed herein is the estimation of systematic model errors in output quantities of interest at future times, and the use of this information to improve the model forecasts. The second problem considered is the identification of those specific physical processes that contribute most to the forecast uncertainty in the quantity of interest under specified meteorological conditions. 
In order to address these questions we employ two machine learning approaches, random forests and artificial neural networks. The discrepancies between model results and observations at past times are used to learn the relationships between the choice of physical processes and the resulting forecast errors.

Numerical experiments are carried out with the Weather Research and Forecasting (WRF) model. The output quantity of interest is the model precipitation, a variable that is both extremely important and very challenging to forecast. The physical processes under consideration include various micro-physics schemes, cumulus parameterizations, short wave, and long wave radiation schemes. The experiments demonstrate the strong potential of machine learning approaches to aid the study of model errors. 

\end{abstract}

\end{frontmatter}
%

%
%

\section{Introduction}
\label{Intro_wrf}

Computer simulation models of the physical world, such as numerical weather prediction (NWP) models, are imperfect and can only approximate the complex evolution of physical reality. Some of the errors are due to the uncertainty in the initial and boundary conditions, forcings, and model parameter values. Other errors, called structural model errors, are due to our incomplete knowledge about the true physical processes, and manifest themselves as missing dynamics in the model \cite{moosavi2017structual}. Examples of structural errors include the misrepresentation of sea-ice in the spring and fall, errors affecting the stratosphere above polar regions in winter \cite{tremolet2007model}, as well as errors due to the interactions among (approximately-represented) physical processes. 

Data assimilation improves model forecasts by fusing information from both model outputs and observations of the physical world in a coherent statistical estimation framework \cite{akella2009different, le1986variational, navon1992variational, tremolet2007model}. While traditional data assimilation reduces the uncertainty in the model state and model parameter values, no methodologies to reduce the structural model uncertainty are available to date.

%

In this study we consider the Weather Research and Forecasting (WRF) model \cite{wrf_website},
a mesoscale atmospheric modeling system. The WRF model includes multiple physical processes and parametrization schemes,
and choosing different model options can lead to significant variability in the model predictions \cite{fovell2006impact, nasrollahi2012assessing}. 

Among different atmospheric phenomena, the prediction of precipitation is extremely challenging and is obtained by solving the atmospheric dynamic and  thermodynamic equations \cite{nasrollahi2012assessing}. Model forecasts of precipitation are very sensitive to physics options such as the micro-physics, cumulus, long wave, and short wave radiation \cite{fovell2010influence,nasrollahi2012assessing, lowrey2008assessing}. Other physics settings that can affect the WRF precipitation predictions include surface physics, 
planetary boundary layer (PBL),  land-surface (LS) parameterizations, and lateral boundary condition. Selecting the right physical process representations and parameterizations is a challenge. In practice the values of physical parameters are empirically determined such as to minimize the difference between the measurements and model predictions \cite{wrf_website, lowrey2008assessing}. 

Considerable effort has been dedicated to determining the best physical configurations of the weather forecast models such as to improve 
their predictions of precipitation. No single choice of physical parameters works perfectly for all times, geographical locations, or meteorological conditions \cite {gallus1999eta, wang1997comparison}. Lowrey and Yang \cite{lowrey2008assessing} investigated the errors in precipitation predictions caused by different parameters including micro-physics and cumulus physics, the buffer zone, the initialization interval, the domain size and the initial and boundary conditions. Jankov et al. \cite{jankov2007impact} examined different combinations of cumulus convection schemes, micro-physical options, and boundary conditions. They concluded that  no configuration was the clear winner at all times, and the variability of precipitation predictions was more sensitive to the choice of the cumulus options rather than micro-physical schemes.
Another study conducted by Nasrollahi \cite{nasrollahi2012assessing} showed that the best model ability to predict hurricanes was achieved using a particular cumulus parameterization scheme combined with a particular micro-physics scheme. Therefore, the interactions of different physical parameterizations have a considerable impact on model errors, and can be considered as one of the main sources of uncertainty that affect the forecast accuracy.

This paper demonstrates the potential of machine learning techniques to help solve two important problems related to the structural/physical uncertainty in numerical weather prediction models. he first problem addressed herein is the estimation of systematic model errors in output quantities of interest at future times, and the use of this information to improve the model forecasts. The second problem considered is the identification of those specific physical processes that contribute most to the forecast uncertainty in the quantity of interest under specified meteorological conditions. 

The application of machine learning techniques to problems in environmental science has grown considerably in recent years. In \cite{gilbert2010machine} a kernel based regression method is developed as a forecasting approach with performance close to Ensemble Kalman Filter (EnKF) and less computational resources. Krasnopol et al. \citep{krasnopol_sky2011development} employ an Artificial Neural Network technique for developing an ensemble stochastic convection parameterization for climate models. Attia et al. \cite{attia2016clusterHMC} develop a new filtering algorithm called Cluster Hybrid Monte Carlo sampling filter (CLHMC) non-Gaussian data assimilation which relaxes the Gaussian assumptions by employing a clustering step.  Moosavi et al. \cite{moosavi2018localization} use regression machine learning techniques for adaptive localization in ensemble based data assimilation. 

This study focuses on the uncertainty in forecasts of cumulative precipitation caused by imperfect representations of physics
and their interaction in the WRF model. The total accumulated precipitation includes all phases of convective and non-convective precipitation.
Specifically, we seek to use the discrepancies between WRF forecasts and measured precipitation levels in the past in order to 
estimate in advance the WRF prediction uncertainty. The model-observation differences contain valuable information about the error dynamics and the missing physics of the model.
We use this information to construct two probabilistic functions. The first one maps the discrepancy data and the physical parameters
onto the expected forecast errors. The second maps the forecast error levels onto the set of physical parameters that are consistent with them. Both maps are constructed using supervised machine learning techniques, specifically, using Artificial Neural Networks and Random Forests  \cite{murphy2012machine}. The two probabilistic maps are used to address the problems posed above, namely the estimation of model errors in output quantities of interest at future times, and the identification of physical processes that contribute most to the forecast uncertainty. 

The remainder of this study is organized as follows. Section \ref{model_error} covers the definition of the model errors.
Section \ref{ML_error_modeling} describes the proposed approach of error modeling using machine learning.
Section \ref{numerical_experiment_wrf} reports numerical experiments with the WRF model that illustrate the capability of the new approach to answer two important questions regarding model errors.
Conclusions are drawn in Section \ref{Conc_wrf}.




\section{Model errors}
\label{model_error}

First-principles computer models capture our knowledge about the physical laws that govern the evolution of a real physical system. The model evolves an initial state at the initial time to states at future times. All models are imperfect, e.g., atmospheric model uncertainties are associated with sub-grid modeling, boundary conditions, and forcings. All these modeling uncertainties are
aggregated into a component that is generically called {\em model
error} \cite{Glimm_2004,Orrell_2001,Palmer_2005}. In the past decade there has been a considerable scientific effort to incorporate model errors and estimate their impact on the best estimate in both variational and statistical approaches
\cite{akella2009different, cardinali2014representing, hansen2002accounting, rao2015posteriori, tr2006accounting, tremolet2007model, zupanski2006model}.

In what follows, we describe our mathematical formulation of the model error associated with NWP models. A similar formulation has been used in \cite{moosavi2017structual} where the model structural uncertainty is studied based on the information provided by the discrepancy between the model solution and the true state of the physical system, as measured by the available observations. 

Consider the following NWP computer model $\mathcal{M}$, that describes the time-evolution of the state of the atmosphere:
\begin{subequations}
\label{eq:Model}
\begin{equation} 
\label{eq:ModelDynamics}
\mathbf{x}_{t} = \mathcal{M}\left(\mathbf{x}_{t-1},  \Uptheta \right), \quad t=1, \cdots, T\,.
\end{equation}
The state vector $\mathbf{x}_t \in \Re^n$ contains the dynamic variables of the atmosphere such as temperature, pressure, precipitation, tracer concentrations, at all spatial locations covered by the model, and at $t$. All the physical parameters of the model are lumped into $\Uptheta \in \Re^\ell$. 

Formally, the true state of the atmosphere can be described by a physical process $\mathcal{P}$ with internal states $\mathbf{\upsilon}_{t}$, which are unknown.  The atmosphere, as an abstract physical process, evolves in time as follows:
\begin{equation} 
\label{eq:NatureDynamics-ideal}
\mathbf{\upsilon}_{t} = \mathcal{P}\left(\mathbf{\upsilon}_{t-1}\right), \quad t=1, \cdots, T.
\end{equation}
The model state seeks to approximates the physical state:
\begin{equation} 
\label{eq:ModelState}
\mathbf{x}_{t} \approx \psi(\mathbf{\upsilon}_t) , \quad t=1, \cdots, T,
\end{equation}
\end{subequations}
where the operator $\psi$ maps the physical space onto the model space, e.g., by sampling the continuous meteorological fields onto a finite dimensional computational grid \cite{moosavi2017structual}. 

Assume that the model state at $t-1$ has the ideal value obtained from the true state via \eqref{eq:ModelState}. The model prediction at $t$ will differ from the reality:
\begin{equation} 
\label{eq:StateDiscrepancy}
 \psi(\mathbf{\upsilon}_{t}) = \mathcal{M}\bigl(\psi(\mathbf{\upsilon}_{t-1}), \Uptheta \bigr)
  + \bdelta_{t}\bigl(\mathbf{\upsilon}_{t}\bigr), \quad t=1, \cdots, T,
\end{equation}
where the discrepancy $\bdelta_t \in \Re^n$ between the model prediction and reality is the {\it structural model error.} This vector lives in the model space. 

Although the global physical state $\mathbf{\upsilon}_{t}$ is unknown, we obtain information about it by measuring of a finite number of observables $\mathbf{y}_t \in \Re^m$, as follows:
\begin{equation} 
\label{eq:NatureObservations-ideal}
\mathbf{y}_t = h(\mathbf{\upsilon}_t) + \epsilon_t,  \quad \epsilon_t \sim \mathcal{N}(0, \mathbf{R}_t), \quad t=1, \cdots, T,
\end{equation}
Here $h$ is the observation operator that maps the true state of atmosphere to the observation space, and the observation error $\epsilon_t$ is assumed to be normally distributed.

In order to relate the model state to observations we also consider the observation operator $\mathcal{H}$ that maps the model state onto the observation space; the model-predicted values $\mathbf{o}_t \in \Re^m$ of the observations \eqref{eq:NatureObservations-ideal} are:
\begin{equation} 
\label{eq:ModelObservations}
\mathbf{o}_t = \mathcal{H}(\mathbf{x}_t), \quad t=1, \cdots, T.
\end{equation}
We note that the measurements $\mathbf{y}_t$ and the predictions $\mathbf{o}_t$ live in the same space and therefore can be directly compared. The difference between the observations \eqref{eq:NatureObservations} of the real system and the model predicted values of these observables \eqref{eq:ModelObservations} represent the model error in observation space: 
\begin{equation} 
\label{eq:ObservedDiscrepancy}
 \boldsymbol{\Delta}_{t} = \mathbf{o}_{t} - \mathbf{y}_{t} \in \Re^m, \quad t=1, \cdots, T.
\end{equation}

For clarity, in what follows we make the following simplifying assumptions \cite{moosavi2017structual}: 
\begin{itemize}
\item the physical system is finite dimensional $\mathbf{\upsilon}_t \in \Re^n$, 
\item the model state lives in the same space as reality, i.e., $\mathbf{x}_{t} \approx \mathbf{\upsilon}_t$ and $\psi(\cdot)\equiv id$ is the identity operator in \eqref{eq:ModelState}, and \item $\mathcal{H}(\cdot) \equiv h(\cdot)$ in \eqref{eq:NatureObservations-ideal} and \eqref{eq:ModelObservations}. 
\end{itemize}
{\it These assumptions imply that the discretization errors are very small, and that  the main source of error are the parameterized physical processes represented by $\Uptheta$ and the interaction among these processes. Uncertainties from other sources, such as boundary conditions, are assumed to be negligible.}

With these assumptions, the evolution equations for the physical system \eqref{eq:NatureDynamics-ideal} and the physical observations equation \eqref{eq:NatureObservations-ideal} become, respectively:
\begin{subequations}
\label{eq:Nature}
\begin{eqnarray} 
\label{eq:NatureDynamics}
\mathbf{\upsilon}_{t} &=& \mathcal{M}\bigl(\mathbf{\upsilon}_{t-1}, \Uptheta \bigr)+ \bdelta_{t}\bigl(\mathbf{\upsilon}_{t}), \quad t=1, \cdots, T, \\
\label{eq:NatureObservations}
\mathbf{y}_t &=& h(\mathbf{\upsilon}_t) + \epsilon_t.
\end{eqnarray}
\end{subequations}

The model errors $\bdelta_t$ \eqref{eq:StateDiscrepancy} are not fully known at any time $t$, as having the exact errors is akin to having a perfect model. However, the discrepancies between the modeled and measured observable quantities \eqref{eq:ObservedDiscrepancy} at past times have been computed and are available at the current time $t$.

Our goal is to use the errors in observable quantities at past times, $\boldsymbol{\Delta}_\tau$ for $\tau=t-1,t-2,\cdots$, in order to estimate the model error $\bdelta_\tau$ at future times $\tau=t,t+1,\cdots$. This is achieved by unravelling the hidden information in the past $\boldsymbol{\Delta}_\tau$ values. Good estimates of the discrepancy $\bdelta_t$, when available, could improve model predictions by applying the correction \eqref{eq:NatureDynamics} to model results: 
\begin{equation}
\label{eq:ModelCorrection}
\mathbf{v}_{t} \approx \mathbf{x}_{t} + \bdelta_{t}.
\end{equation} 

Our proposed error modeling approach constructs input-output mappings to estimate given aspects of model errors $\bdelta_t$. The inputs to these mappings are the physical parameters $\Uptheta$ of the model. The outputs to these mappings are different aspects of the error in a quantity of interest, such as the model errors over a specific geographical location, or the error norm of model error integrated over the entire domain.

Specifically, the aspect of interest (quantity of interest) in this study is the error in precipitation levels forecasted by the model. The parameters $\Uptheta$ describe the set of physical processes that are essential to be included in the WRF model in order to produce accurate precipitation forecasts. The WRF model is modular and different combinations of the physical packages can be selected, each corresponding to a different value of $\Uptheta$.

We use the error mappings learned from past model runs to estimate the model errors $\bdelta_{t}$.
We also consider estimating what combination of physical processes $\Uptheta$ leads to lower model errors, or reversely, what interactions of which physics cause larger errors in the prediction of the quantity of interest. 

\section{Approximating model errors using machine learning}
\label{ML_error_modeling}
We propose a multivariate input-output learning model to predict the model errors $\bdelta$, defined in \eqref{eq:StateDiscrepancy}, stemming from the uncertainty in parameters $\Uptheta$.
To this end, we define a probabilistic function $\phi$ that maps every set of input features $F \in \Re^r$ to output target variables $\Uplambda \in \Re^o$: 
\begin{equation}
\label{eqn:general_learning}
\phi ( F ) \approx \Uplambda\,,
\end{equation}
and approximate the function $\phi$ using machine learning.

Different particular definitions of $\phi$ in \eqref{eqn:general_learning} will be used to address two different problems related to model errors, as follows:
\begin{enumerate}
\item The first problem is to estimate the systematic model error in certain quantities of interest at future times, and to use this information in order to improve the WRF forecast. To achieve this one quantifies the model error aspects that correspond to running WRF with different physical configurations (different parameters $\Uptheta$).
\item The second problem is to identify the specific physical processes that contribute most to the forecast uncertainty in the quantity of interest under specified meteorological conditions. To achieve this one finds the model configurations (physical parameters $\Uptheta$) that lead to forecast errors smaller that a given threshold under specified meteorological conditions.
\end{enumerate}
In what follows we explain in detail the function $\phi$ specification, the input features, and the target variables for each of these problems.
\subsection{Problem one: estimating in advance aspects of interest of the model error} 
\label{sect:error_model_predict}
Forecasts produced by NWP models are contaminated by model errors. These model errors are highly correlated in time; hence historical information about the model errors can be used as an input to the learning model to gain insight about model errors that affect the forecast. We are interested in the uncertainty caused due to the interaction between the various components in the physics based model; these interactions are lumped into the parameter $\Uptheta$ that is supplied as an input to the learning model. The learning model aims to predict the error of NWP model of next forecast window using the historical values of model error and the physical parameters used in the model. We define the following mapping:
 \begin{equation}\label{eqn:error_mapping}
 \phi^{\rm error} \left( \Uptheta, \boldsymbol{\Delta}_{\tau}, \mathbf{o}_{\tau}, \mathbf{o}_{t} \right) \approx \boldsymbol{\Delta}_{t}\,
 \quad \tau< t.
\end{equation}
We use a machine learning algorithm to approximate the function $\phi^{\rm error}$. 
The learning model is trained using a dataset that consists of the following inputs:
\begin{itemize}
        \item WRF physical packages that affect the physical quantity of interest ($\Uptheta$),
        \item historical WRF forecasts ($\mathbf{o}_{\tau}$ for $\tau \le t-1$), 
        \item historical model discrepancies ($\mathbf{\Delta}_{\tau}$ for $\tau \le t-1$), 
        \item WRF forecast at the current time ($\mathbf{o}_t$), 
        \item the available model discrepancy at the current time ($\boldsymbol{\Delta}_t$) since we have access to the 
         observations from reality $y_t$ at the current time step.
\end{itemize}
In supervised learning process, the learning model identifies the effect of physical packages, the historical WRF forecast, the historical model discrepancy, and the WRF forecast at the current time on the available model discrepancy at the current time.
After the model get trained on the historical data, it yields an approximation to the mapping $\phi^{\rm error}$. 
We denote this approximate mapping by $\widehat{\phi}^{\rm error}$.

During the test phase the approximate mapping $\widehat{\phi}^{\rm error}$ is used to estimate the model discrepancy $\widehat{\boldsymbol{\Delta}}_{t+1}$ in advance. We emphasize that the model prediction (WRF forecast) 
at the time of interest $t+1$ ($\mathbf{o}_{t+1}$) is available, where as the model discrepancy $\widehat{\boldsymbol{\Delta}}_{t+1}$ is an unknown quantity.
In fact the run time of WRF is much smaller than the time interval
between $t$ and $t+1$, or in other way, the time interval is large enough to run the WRF model and obtain the forecast for next time window, estimate the model errors for next time window and finally improve the model forecast by combining the model forecast and 
model errors. 

At the test time we predict the future model error as follows:
 \[
 \widehat{\boldsymbol{\Delta}}_{t+1} \approx \widehat{\phi}^{\rm error} \left( \Uptheta,\boldsymbol{\Delta}_{\tau}, \mathbf{o}_{\tau}, \mathbf{o}_{t+1}  \right)\,, \quad \tau<t+1\,.
\]

As explained in \cite{moosavi2017structual}, the predicted error $\widehat{\boldsymbol{\Delta}}_{t+1}$ in the observation space can be used to estimate the error $\bdelta_{t+1}$ in the model space. In order to achieve this one needs to use additional information about the structure of the model and the observation operator.
For example, if the error $\widehat{\boldsymbol{\Delta}}_{t+1}$ represents the projection of the full model error onto the observation space, we have:
\begin{subequations}
\label{eqn:error-obs-to-model}
\begin{equation}
\label{eqn:error-obs-to-model-H}
\boldsymbol{\Delta}_{t+1} \approx \mathbf{H}_t \cdot\bdelta_{t+1}, \quad \widehat{\bdelta}_{t+1} \approx \mathbf{H}_t\,\left( \mathbf{H}_t^T\, \mathbf{H}_t\right)^{-1}\,\mathbf{H}_t^T \cdot \widehat{\boldsymbol{\Delta}}_{t+1},
\end{equation}
where we use the linearized observation operator at the current time, $\mathbf{H}_t = h'(\boldsymbol{x}_{t})$. A more complex approach is to use a Kalman update formula:
\begin{equation}
\label{eqn:error-obs-to-model-K}
\widehat{\bdelta}_{t+1} \approx \textnormal{cov}(\mathbf{x}_t,\mathbf{o}_t)\, \left(\textnormal{cov}(\mathbf{o}_t,\mathbf{o}_t)  + \mathbf{R}_t\right)^{-1}\, \widehat{\boldsymbol{\Delta}}_{t+1},
\end{equation}
\end{subequations}
where $\mathbf{R}_t$ is the covariance of observation errors. The Kalman update approach requires estimates of the covariance matrices between model variables; such covariances are already available in an ensemble based data assimilation system. Once we estimate the future model error $\bdelta_{t+1}$, we can improve the NWP output using equation \eqref{eq:ModelCorrection}.

\subsection{Problem two: identifying the physical packages that contribute most to the forecast uncertainty}
Typical NWP models incorporate an array of different physical packages to represent multiple physical phenomena that act simultaneously. Each physical package contains several alternative configurations (e.g., parameterizations or numerical solvers) that affect the accuracy of the forecasts produced by the NWP model. A particular scheme in a certain physical package best captures the reality under some specific conditions (e.g., time of the year, representation of sea-ice, etc.). The primary focus of this study is the accuracy of precipitation forecasts, therefore we seek to learn the impacts of all the physical packages that affect precipitation. To this end, we define the following mapping:

 \begin{equation}\label{eqn:physics_mapping}
    \phi^{\rm physics}\left( \boldsymbol{\Delta}_{t} \right) \approx \Uptheta\,,
\end{equation}
that estimates the configuration $\Uptheta$ of the physical packages such that the WRF run generates a forecast with an error consistent with the prescribed level $\boldsymbol{\Delta}_{t}$ (where $\boldsymbol{\Delta}_{t}$ defined in equation \eqref{eq:ObservedDiscrepancy} is the forecast error in observation space at time $t$.)

We train the model to learn the effect of the physical schemes on the mismatch between WRF forecasts and reality. The input data required for the training process is obtained by running the model with various physical package configurations $\Uptheta^{\rm train}_i$, and comparing the model forecast against the observations at all past times $\tau$ to obtain the corresponding errors $\boldsymbol{\Delta}_{\tau,i}^{\rm train}$ for $\tau \le t$ and $i \in \{training~data~set\}$. The output data is the corresponding 
physical combinations $\Uptheta$ that leads to the input error threshold.

In order to estimate the combinations of physical process configuration that contribute most to the uncertainty in predicting precipitation we take the following approach.  The dataset consisting of the  observable discrepancies during the current time window $\boldsymbol{\Delta}_{t}$ 
is split into a training part and a testing part.  In the test phase we use the approximated function $\widehat{\phi}^{\rm physics}$ to estimate the physical process settings $\widehat{\Uptheta}_j^1$ that are consistent with the observable errors $\boldsymbol{\Delta}_{t,j}^{\{1\}}$. Here we select  $\boldsymbol{\Delta}_{t,j}^{\{1\}} = \boldsymbol{\Delta}_{t,j}^{\rm test}$ for each $j \in \{test~data~set\}$. Note that in this case, since we know what physics has been used for the current results,  one can take $\widehat{\Uptheta}_j^{\{1\}}$ to be the real parameter values $\Uptheta_j^{\{1\}}$ used to generate the test data. However, in general, one selects $\boldsymbol{\Delta}_{t,j}^{\{1\}}$ in an application-specific way and the corresponding parameters need to be estimated.

Next, we reduce the desired forecast error level to $\boldsymbol{\Delta}_{t,j}^{\{2\}}=\boldsymbol{\Delta}_{t,j}^{\{1\}}/2$, and use the approximated function $\widehat{\phi}^{\rm physics}$ to estimate the physical process setting $\widehat{\Uptheta}_j^{\{2\}}$ that corresponds to this more accurate forecast. To identify the package setting that has the largest impact on the observable error we monitor the variability in the predicted parameters $\widehat{\Uptheta}^{\{2\}} - \widehat{\Uptheta}^{\{1\}}$. Specifically,  the number of times the setting of a physical process in $\widehat{\Uptheta}_j^2$ is different from its setting in $\widehat{\Uptheta}_j^1$ is an indicator of the variability in model prediction when that package is changed. A higher variability in predicted physical packages implies a larger contribution towards the model errors - as estimated by the ML model.

\subsection{Machine learning algorithms}
\label{Machine_learning}

In order to approximate the functions \eqref{eqn:error_mapping} and \eqref{eqn:physics_mapping} discussed earlier we use regression machine learning methods. Choosing a right learning algorithm to use is challenging as it largely depends on the problem and the data available~\cite{attia2016clusterHMC,asgari2017Latent,moosavi2017multivariate,moosavi2018localization}. Here, we use Random Forests (RF) and Artificial Neural Networks (ANN) as our learning algorithms \cite{murphy2012machine}. Both RF and ANN algorithms tan handle non-linearity in regression and classification. Given that the physical phenomena governing precipitation are highly nonlinear, and and atmospheric dynamics is chaotic, we believe that RF and ANN approaches are well suited to capture the associated features.
We briefly review these techniques  next.

\subsubsection{Random forests}
\label{random_2orest}
A random forest \cite{breiman2001random} is an ensemble based method that constructs multiple decision trees. The  principle idea behind ensemble methods is that a group of weak learners can come together to form a strong learner \cite{ breiman1996bagging, breiman2001random}.  The decision tree is built top-down from observations of target variables. The observation dataset is partitioned, smaller subsets are represented in branches, and decisions about the target variables are represented in the leaves.

There are many specific decision-tree algorithms available, including ID3 (Iterative Dichotomiser 3) 
\cite{quinlan1986induction}, C4.5 (successor of ID3) \cite{quinlan2014c4}, 
CART (Classification And Regression Tree), 
CHAID (CHi-squared Automatic Interaction Detector), 
and conditional inference trees~\cite{strobl2008conditional}.
If the dataset has multiple attributes, one can decide which attribute to place at the root or at different levels of the tree 
by considering different criteria such as information gain or the gini index \cite{ceriani2012origins}. 

Trees can be non-robust, with small changes in the tree leading to large changes in regression results. Moreover, trees tend to over-fit the data \cite{segal2004machine}. The random forest algorithm uses the bagging technique for building an ensemble
of decision trees which are accurate and powerful at handling large, high dimensional datasets. Moreover, the bagging 
technique greatly reduces the variance ~\cite{dietterich2000ensemble}.
For each tree in the forest, a bootstrap sample \cite{breiman1996bagging, dietterich2000ensemble} is selected from the dataset and instead of examining all possible feature-splits, some subset of the features is selected \cite{liaw2002classification}.
The node then splits on the best feature in the subset.  By using a random sample of features the correlation between trees in the ensemble
decreases, and the learning for each tree  is much faster by restricting the features considered for each node. 

\subsubsection{Artificial neural networks}
\label{sect:NN}
ANN is a computational model inspired by human brain's biological structure.  
ANN consist of neurons and connections between the neurons (weights) which are organized in layers.
At least three layers of neurons (an input layer, a hidden layer, and an output layer) are required for construction of a neural network, where
the input layer distributes the input signals to the first hidden layer. The feed-forward operation in a network passes information to neurons in a subsequent hidden layer. The neurons combine this information, and the output of each layer is obtained by passing the combined information through  a differentiable transfer function that can be log-sigmoid, hyperbolic tangent sigmoid, or linear transfer function.

In supervised learning the network is provided with samples from which it discovers the relations of inputs and outputs. 
The learning problem consists of finding the optimal parameters of network such that 
the error between the desired output and the output signal of the network is minimized.
The network first is initialized with randomly chosen weights and then 
the error is back-propagated through the network using a gradient descent method. 
The gradient of the error function is computed and used to modify weights and biases such that the error between the desired output and the output signal of the network is minimized \cite{funahashi1989approximate, rumelhart1985learning} .
This process is repeated iteratively until the network output is close to the desired output \cite{haykin2009neural}.

\section{Numerical experiments}
\label{numerical_experiment_wrf}
We apply the proposed learning models to the Weather Research and Forecasting model \cite{wrf_website} in order to: 
\begin{itemize}
  \item predict the bias in precipitation forecast caused by structural model errors,
  \item predict the statistics associated with the precipitation errors, and
  \item identify the specific physics packages that contribute most to precipitation forecast errors for given meteorological conditions.
\end{itemize}

\subsection{The WRF model}
\label{wrf}

\begin{figure}[H]
\centering
\subfigure[NCEP analysis at $6 \textsc{am}$ provides initial conditions]{
\includegraphics[scale=0.4]{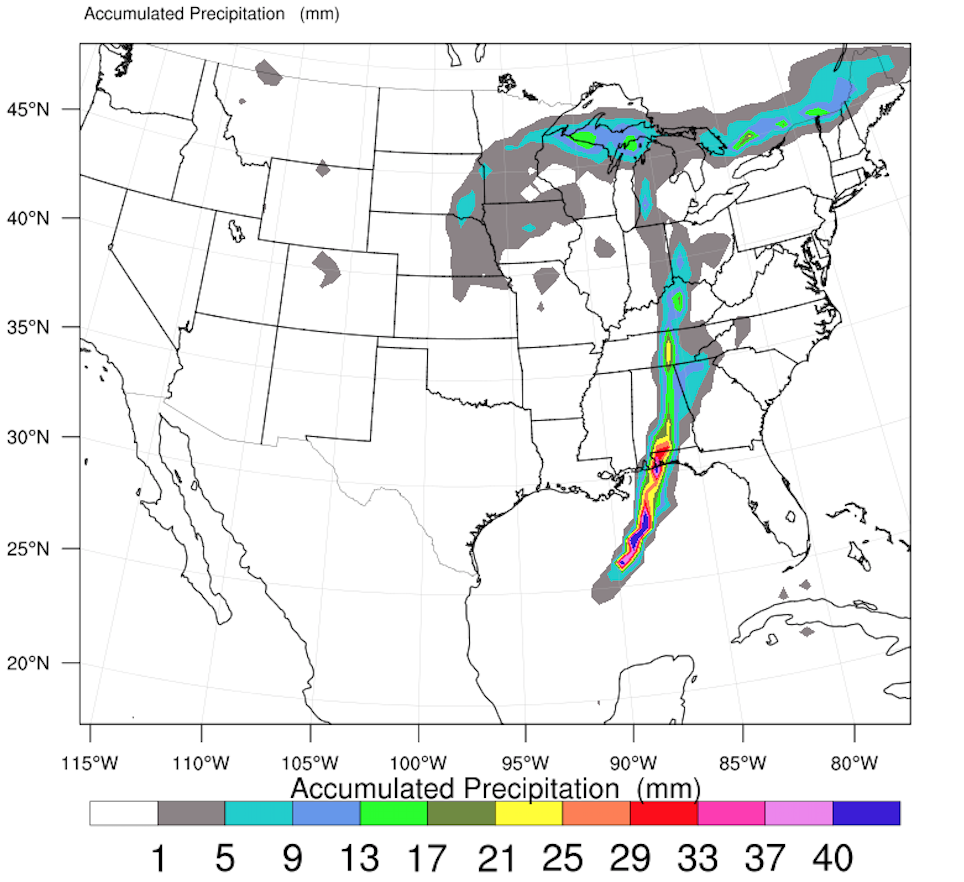}
\label{fig:NCEP_initial}
}
\subfigure[NCEP analysis at $12 \textsc{pm}$ provides a proxy for the true state of the atmosphere]{
\includegraphics[scale=0.4]{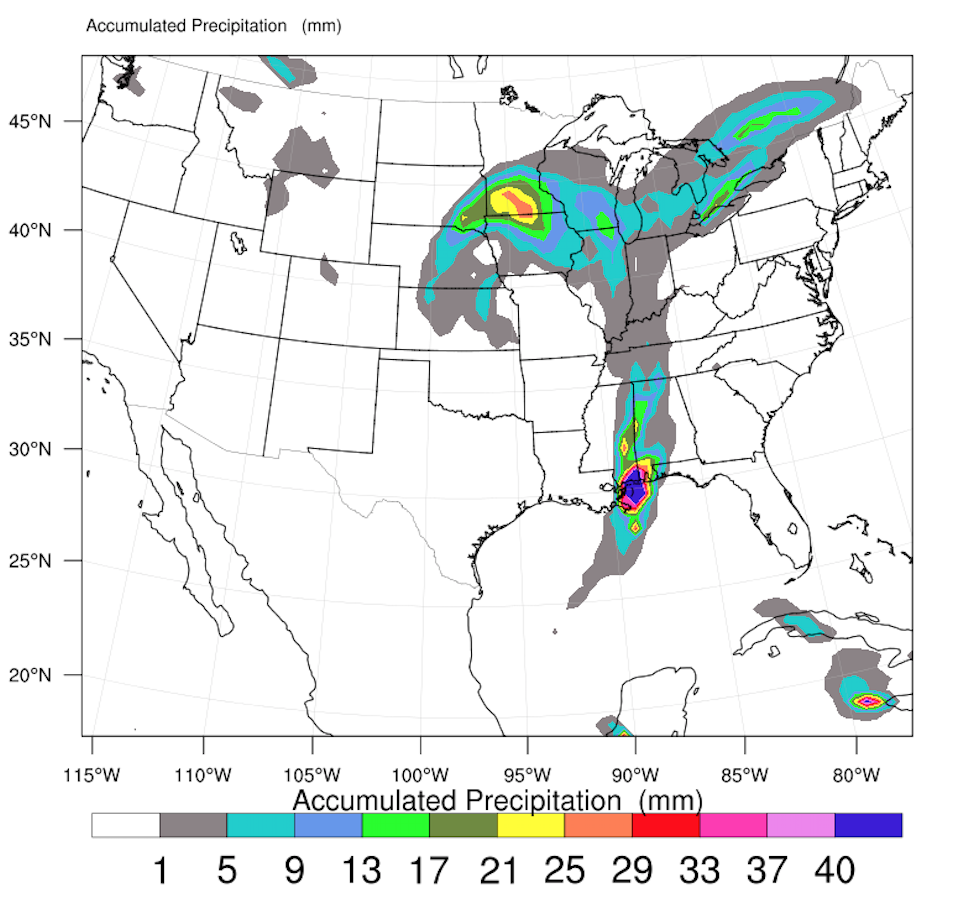}
\label{fig:NCEP_analysis}
}
\subfigure[WRF forecast at $12 \textsc{pm}$ corresponding to the physics micro-physics: Kessler, cumulus physics: Kain-Fritsch,  long wave radiation physics: Cam, short wave radiation physics: Dudhia]{
\includegraphics[scale=0.4]{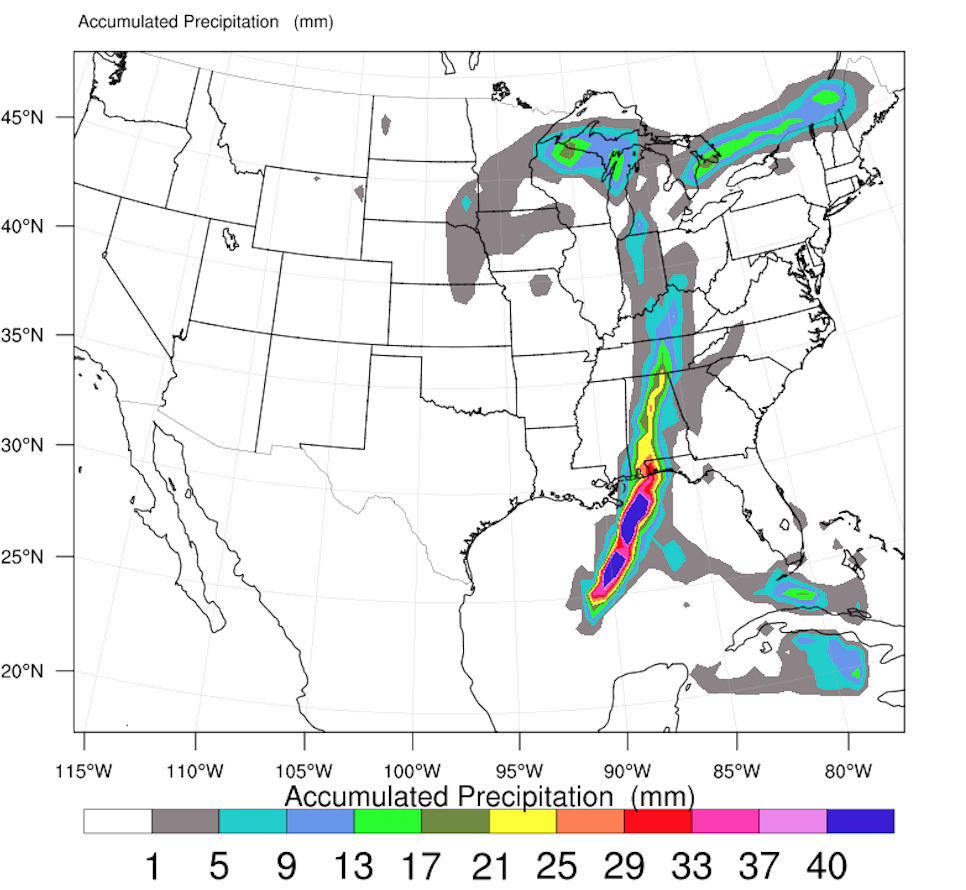}
\label{fig:WRF_output}
}
        \caption{Initial conditions, the analysis and the WRF forecast for the simulation time $12 \textsc{pm}$ on 5/1/2017. Shown in the plots are the accumulated precipitation in millimeter unit.}
\label{fig:NCEP_plots}
\end{figure}

In this study we use the non-hydrostatic WRF model version 3.3. The simulation domain, shown in Fig. \ref{fig:NCEP_plots}, covers the  continental United States and has dimensions of $60 \times 73$ horizontal grid points in the west-east and south-north directions respectively, with a horizontal grid spacing of $60 km$ \cite{wang2014downscaling}.  The grid has 60 vertical levels to cover the troposphere and lower part of the stratosphere between the surface to approximately $20 km$. 
In all simulations, the 6-hourly analysis from the National Centers for Environmental Prediction (NCEP) are used as the initial and boundary conditions of the model \cite{NCEP_data}.
The stage IV estimates are available at an hourly temporal resolution over continental United States. For experimental purposes, we use the stage IV NCEP analysis 
as a proxy for the true state of the atmosphere.
The simulation window begins at $6 \textsc{am}$ UTC (Universal Time Coordinated)  on May 1st 2017, and
the simulation time is a six hour window time the same day. 
The ``true" states of the atmosphere are available using the NCEP analysis data hourly. All the numerical experiments use the 
NCEP analysis data to run WRF model on  May 1st 2017.


The model configuration parameters $\Uptheta$ represent various combinations of micro-physics schemes, cumulus parameterizations, short wave, and long wave radiation schemes. Specifically, each process is represented by the schema values of each physical parameter it uses, as detailed in WRF model physics options and references \cite{wrf_physics}.
The micro-physics option provides
atmospheric heat and moisture tendencies in atmosphere which also accounts
for the vertical flux of precipitation and the
sedimentation process.
The cumulus parameterization
is used to vertically redistribute heat and moisture independent
of latent heating due to precipitation.
The long wave radiation considers clear-sky and cloud upward and downward
radiation fluxes and the short wave radiation considers clear-sky and cloudy solar fluxes.

A total number of  252 combinations of the four physical modules are used in the simulations.
The micro-physics schemes include:
Kessler \cite{kessler1995continuity}, Lin \cite{lin1983bulk}, WSM3 Hong \cite{hong2004revised}, 
WSM5 Hong \cite{hong2004revised}, Eta (Ferrier), WSM6 \cite{hong2006wrf},
Goddard \cite{tao1989ice}, Thompson \cite{thompson2008explicit}, 
Morrison \cite{morrison2009impact}.
The cumulus physics schemes applied are:  
Kain-Fritsch \cite{kain2004kain}, Betts-Miller-Janjic \cite{janjic1994step},  Grell Freitas\cite{grell2013scale}.  
The long wave radiation physics include:
RRTM \cite{mlawer1997radiative}, Cam \cite{conley2012description}.
Short wave radiation physics include:
Dudhia \cite{dudhia1989numerical}, 
Goddard \cite{chou1999solar}, 
Cam \cite{conley2012description}.

\begin{figure}[H]
\centering
\subfigure[Micro-physics scheme: Kessler, cumulus physics: Kain-Fritsch, short wave radiation: Cam, long wave radiation: Dudhia]{
\includegraphics[scale=0.4]{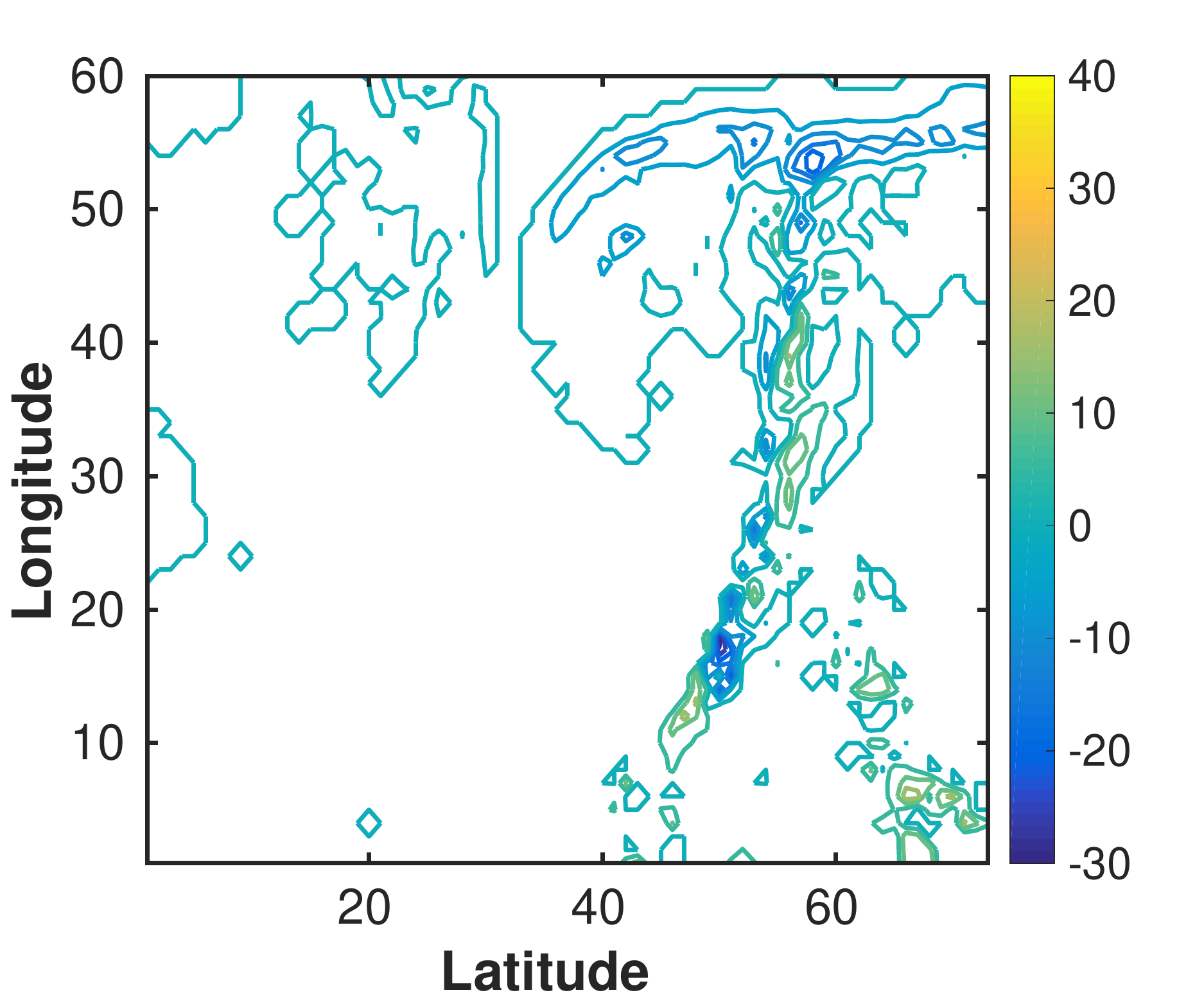}
\label{fig:contour_error_phys1}
}
\subfigure[micro-physics scheme: Lin,  cumulus physics: Kain-Fritsc, short wave radiation: RRTM Mlawer, long wave radiation: Cam]{
\includegraphics[scale=0.4]{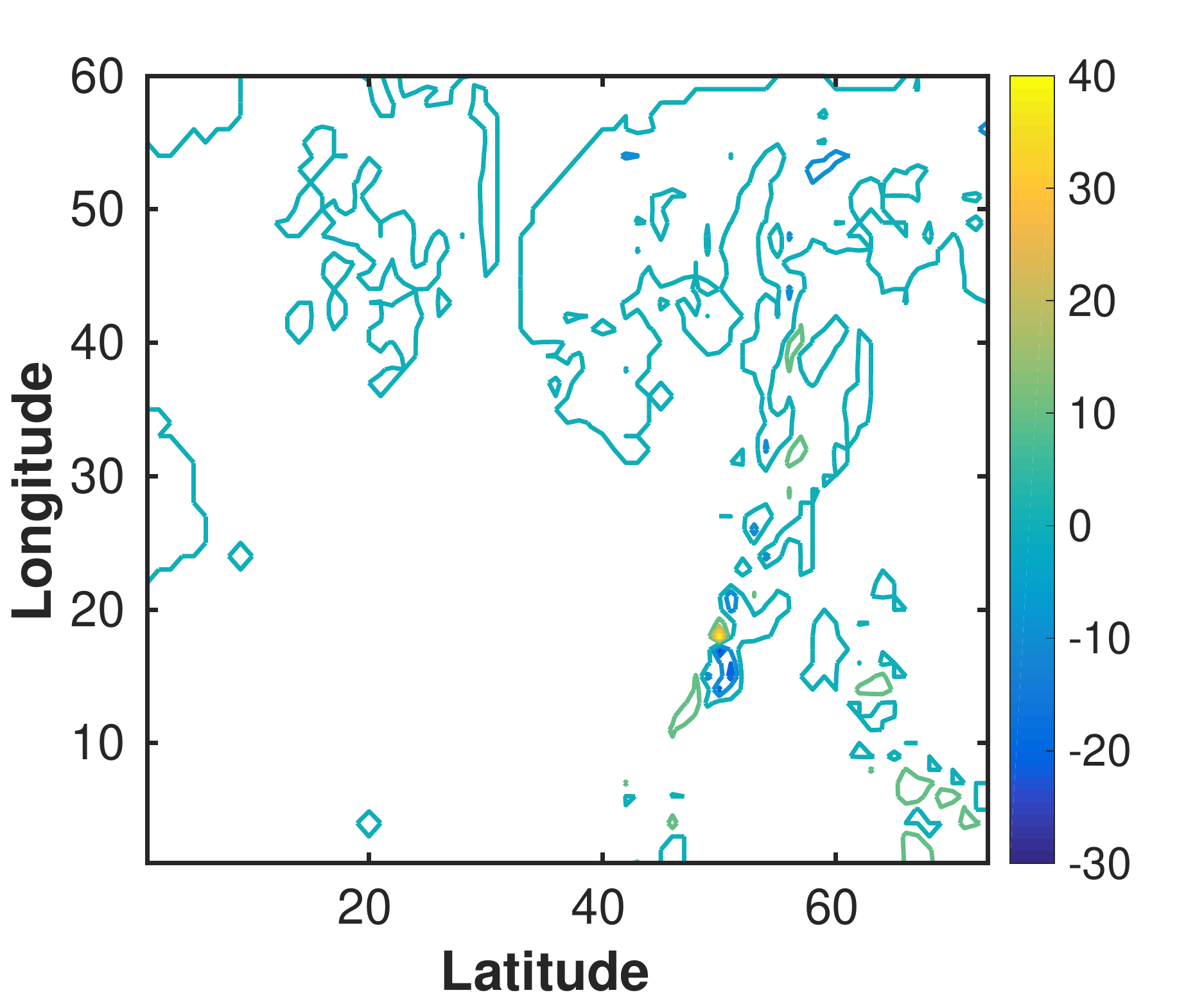}
\label{fig:contour_error_phys2}
}
\caption{
Shown in the plots are contours of 
observable discrepancies which are the differences in the accumulated precipitation results of WRF forecast against the analysis data  $\boldsymbol{\Delta}_{t=12 \textsc{pm}}$ on 5/1/2017 for two different physics combinations. The observation operator extracts the precipitation solution from the WRF state vector.}
\label{fig:contour_error}
\end{figure}

For each of the 252 different physics combinations, the effect of each physics combination on precipitation is investigated.
The NCEP analysis grid points are $428\times 614$, while the WRF computational model have $60 \times 73$ grid points. For obtaining the discrepancy between the WRF forecast and NCEP analysis we linearly interpolate the analysis to transfer the physical variables onto the model grid.  Figure \ref{fig:NCEP_initial} and \ref{fig:NCEP_analysis} shows the NCEP analysis at $6 \textsc{am} $ and $12 \textsc{pm}$ on 5/1/2017 which are used as initial condition and ``true'' (verification) state, respectively. 
The WRF forecast corresponding to the physics micro-physics: Kessler, cu-physics: Kain-Fritsch, ra-lw-physics: Cam , ra-sw-physics: Dudhia is illustrated in Figure \ref{fig:WRF_output}. 
Figure \ref{fig:contour_error} shows contours of discrepancies at $12 \textsc{pm}$ 
$\left( \boldsymbol{\Delta}_{t=12 \textsc{pm}} \right) $ discussed in equation \eqref{eq:ObservedDiscrepancy} for two different physical combinations, which illustrates the effect that changing the physical schemes has on the forecast.

%
%
\subsection{Experiments for problem one: predicting pointwise precipitation forecast errors over a small geographic region}
\label{error_model}
We demonstrate our learning algorithms to forecast precipitation in the state of Virginia on May 1st 2017 at $6 \textsc{pm}$. 
Our goal is to use the learning algorithms to correct the bias created due to model errors and hence improve the forecast for precipitation.
As described in section \ref{sect:error_model_predict}, we learn the function $\phi^{\rm error}$ of equation \eqref{eqn:error_mapping} using the  training data from the previous forecast window ($6 \textsc{am}$ to $12 \textsc{pm}$):
\[
\phi^{\rm error} \left( \Uptheta, \boldsymbol{\Delta}_{\tau}, \mathbf{o}_{\tau}, \mathbf{o}_{t=12 \textsc{pm}} \right) \approx \boldsymbol{\Delta}_{t=12 \textsc{pm}},  \quad   7 \textsc{am} \leq \tau < 12 \textsc{pm}. 
\]
We use two learning algorithms to approximate the function $\phi^{\rm error}$, namely, the RF 
and ANN using Scikit-learn, machine learning library in Python \cite{scikit-learn}.
The RF with ten trees and CART learning tree algorithm is used.
The ANN with six hidden layers and hyperbolic tangent sigmoid activation function in each layer and linear activation function 
at last layer is employed. The number of layers and number of neurons in each layer are tuned empirically.
For training purposes, we use the NCEP analysis of the May 1st 2017 at $6 \textsc{am}$ as initial conditions for the WRF model.
The forecast window is 6 hours and the WRF model forecast final simulation time is $12 \textsc{pm}$. 
The input features are:
\begin{itemize}
\item The physics combinations ($\Uptheta$).
\item The hourly WRF forecasts projected onto observation space $o_{\tau}$, $\textsc{am} \leq \tau \le 12 \textsc{pm}$.
The WRF state ($\mathbf{x}_{t}$) includes all model variables such as temperature, pressure, precipitation, etc.  The observation operator extracts the precipitation portion of the WRF state vector, $\mathbf{o}_{t}\equiv\mathbf{x}_t^\textrm{precipitation}$. Accordingly, $\boldsymbol{\Delta}_{t}$ is the discrepancy between WRF precipitation forecast $\mathbf{o}_{t}$ and the observed precipitation $\mathbf{y}_{t}$.
\item The observed discrepancies at past times ($\boldsymbol{\Delta}_{\tau}$, $7 \textsc{am} \leq \tau < 12 \textsc{pm}$).
\end{itemize}
The output variable is the discrepancy between the NCEP analysis and the WRF forecast at $12 \textsc{pm}$, i.e., the observable discrepancies for the current forecast window ($\boldsymbol{\Delta}_{t=12 \textsc{pm}} $). 
In fact, for each of the $252$ different physical configurations, the WRF model forecast as well as the difference between the WRF forecast and the analysis are provided as input-output combinations for learning the function $\phi^{\rm error}$. The number of grid points over the state of Virginia is $14 \times 12$. Therefore for each physical combination we have $168$ grid points, and the total number of samples in the training data set is $252 \times 168 =42,336$ with $15$ features. 

Both ANN and RF are trained with the above input-output combinations described above
and during the training phase, the learning model learns the 
effect of interaction between different physical configurations on the WRF forecast and model error and obtains the approximation to the function $\phi^{\rm error}$ which we denote by $\widehat{\phi}^{\rm error}$. 
The goal is to have more accurate forecast in the future time windows. We don't have the analysis data of future time windows but we can
run WRF for future time windows and also predict the future model error using the approximated function $\widehat{\phi}^{\rm error}$.
Once we obtain the predicted model error we can use that information in order to raise the accuracy of WRF forecast.
In the testing phase we use the function $\widehat{\phi}^{\rm error}$ to predict the future forecast error $\widehat{\boldsymbol{\Delta}}_{t=6 \textsc{pm}}$ given the combination of physical parameters as well as the WRF forecast at time $6 \textsc{pm}$ as input features. 

\[
\widehat{\boldsymbol{\Delta}}_{t=6 \textsc{pm}} \approx \widehat{\phi}^{\rm error} \left( \Uptheta,\boldsymbol{\Delta}_{\tau}, \mathbf{o}_{\tau}, \mathbf{o}_{t=6 \textsc{pm}}  \right)\,, \quad 1 \textsc{pm} \leq \tau < 6 \textsc{pm}.
\]

To quantify the accuracy of the predicted error we calculate the Root Mean Squared Error (RMSE) between the true and predicted discrepancies at $6\textsc{pm}$:
\begin{equation}\label{eqn:RMSE}
	RMSE = \sqrt{\frac{1}{n} \sum_{i=1}^{n} \left( \widehat{\boldsymbol{\Delta}}_{t=6 \textsc{pm}}^{i}-\boldsymbol{\Delta}_{t=6 \textsc{pm}}^{i}\right)^2 },
\end{equation}
where $n=168$ is the number of grid points over Virginia,. $\widehat{\boldsymbol{\Delta}}_{t=6\textsc{pm}}^{i}$ is the
predicted discrepancy in the $i^{th}$ grid point, and $\boldsymbol{\Delta}_{t=6\textsc{pm}}^{i}$ is the $i^{th}$ actual discrepancy  in the $i^{th}$ grid point. The actual discrepancy is obtained as the difference between the NCEP analysis and WRF forecast at time $t=6\textsc{pm}$. This error metric is computed for each of the $252$ different configurations of the physics. The minimum, maximum and average RMSE over the $252$ runs is reported in Table \ref{tab:RMSE_error}.

\begin{table}[H]
\begin{center}
   \begin{tabular}{| l | l | l | l | l |}
    \hline
   &  $ \textnormal{minimum} (RMSE) $  & $\textnormal{average}(RMSE)$ & $\textnormal{maximum}(RMSE)$
      \\ \hline
ANN   & $ 1.264 \times 10^{-3}$ & $  1.343 \times 10^{-3} $ & $ 5.212 \times 10^{-3}$
     \\ \hline
RF   & $ 1.841 \times 10^{-3}$  & $ 1.931 \times 10^{-3} $ & $ 7.9  \times 10^{-3}$
 \\ \hline
     \end{tabular}
\end{center}
 \caption{The minimum, average, and maximum RMSE between the predicted $\widehat{\boldsymbol{\Delta}}_{t=6 \textsc{pm}}$ and the true $\boldsymbol{\Delta}_{t=6 }$ over 252 physics combinations.}
  \label{tab:RMSE_error}
\end{table}

\begin{figure}[H]
\centering
\subfigure[Original WRF prediction]{
\includegraphics[width=5.7cm, height=6cm]{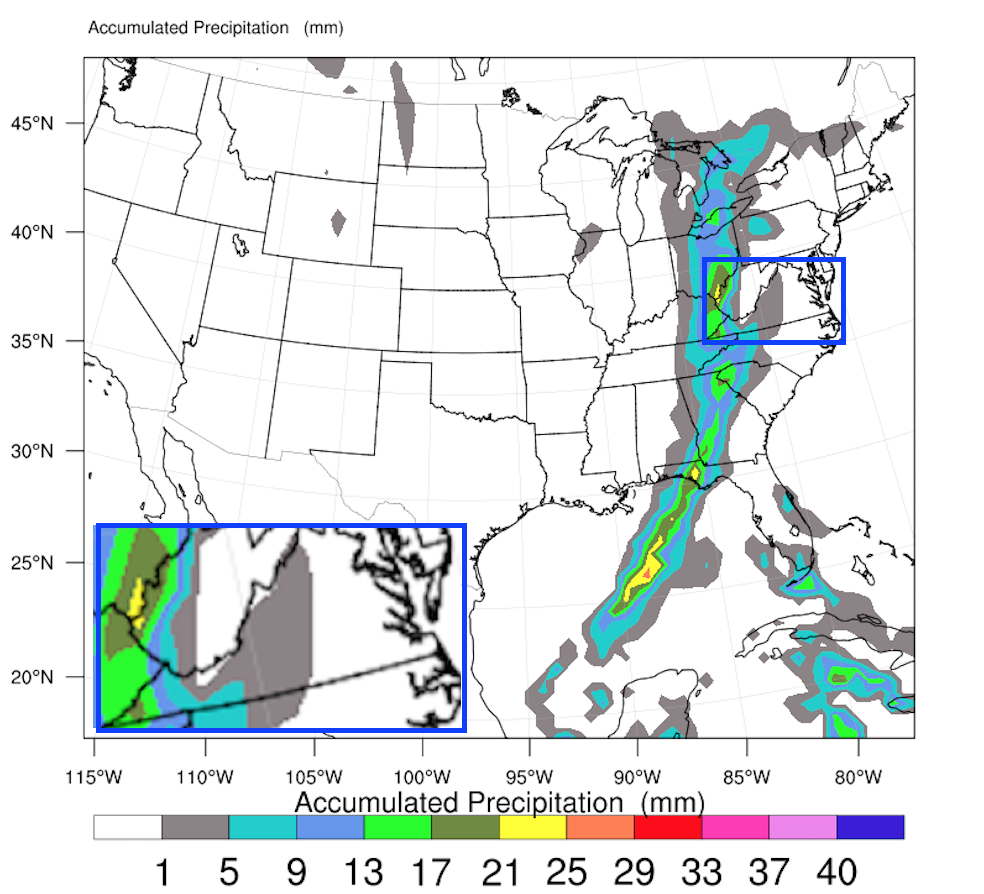}
\label{fig:virginia_phys12}
}
\subfigure[NCEP analysis]{
\includegraphics[width=5.7cm, height=6cm]{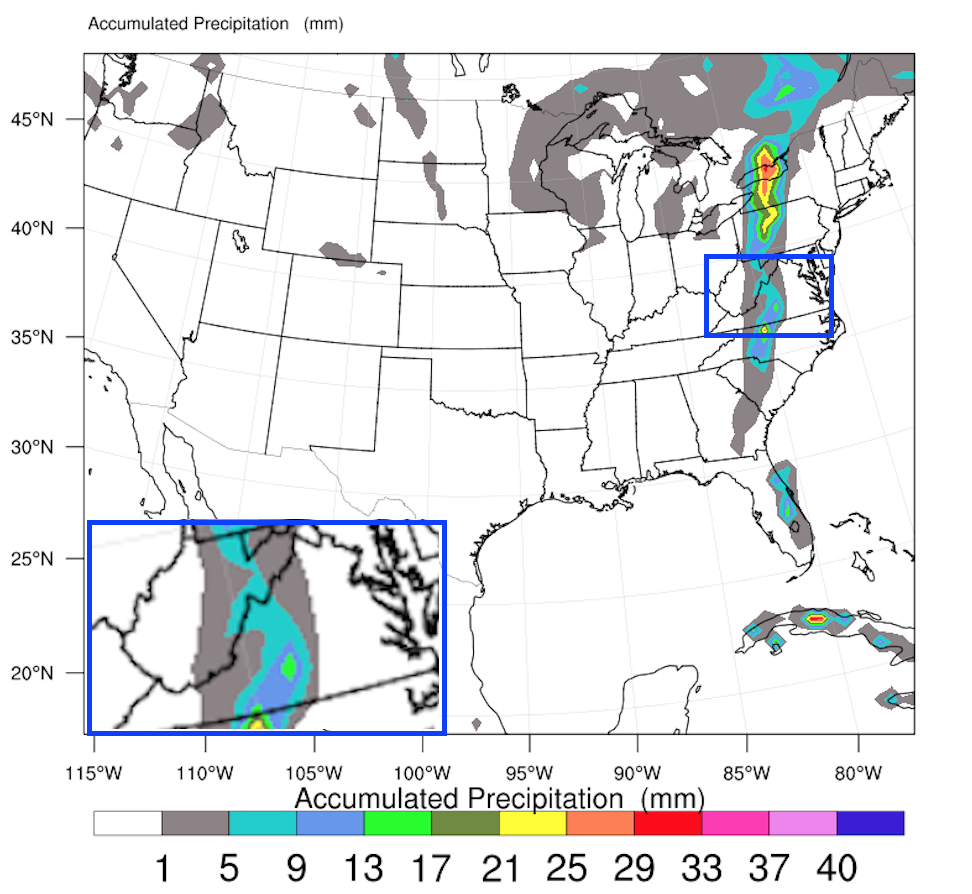}
\label{fig:virginia_analysis}
}
\caption{WRF prediction and NCEP analysis at $6 \textsc{pm}$ on 5/1/2017. Zoom-in panels show the predictions over Virginia.}
\label{fig:virginia_error}
\end{figure}

\begin{figure}[H]
\centering
\subfigure[Discrepancy between original WRF forecast and NCEP analysis]{
\includegraphics[width=5.1cm, height=4.3cm]{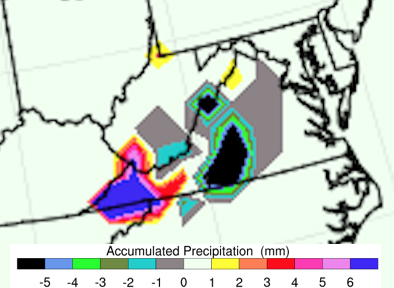}
\label{fig:virginia_analysis_err}
}
\subfigure[Discrepancy between the corrected WRF forecast and the NCEP analysis]{
\includegraphics[width=5.1cm, height=4.3cm]{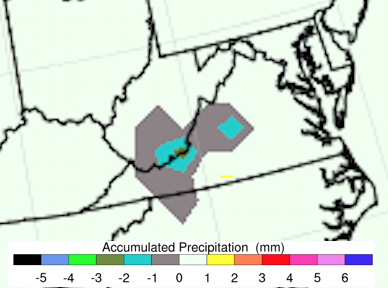}
\label{fig:virginia_improved_err}
}
\caption{Discrepancy between WRF forecasts and the NCEP analysis over Virginia at $6 \textsc{pm}$ on 5/1/2017. The forecast correction clearly improves the model results.}
\label{fig:virginia_improved}
\end{figure}

The predicted discrepancy in the observation space $\widehat{\boldsymbol{\Delta}}_{t=6\textsc{pm}}$  can be used to approximate the discrepancy in the model space $\widehat{\bdelta}_{t=6\textsc{pm}}$ using equation \eqref{eqn:error-obs-to-model}. Here all the grid points are observed and therefore the error in the model space equal to the error in the observation space. Next, the estimate forecast error can be used to correct the forecast bias caused by model errors using \eqref{eq:ModelCorrection}, and hence to improve the forecast at $6\textsc{pm}$:
$\widehat{\mathbf{x}}_{t=6\textsc{pm}} = \mathbf{x}_{t=6\textsc{pm}} + \widehat{\bdelta}_{t=6\textsc{pm}}$.
Figure \ref{fig:virginia_phys12} shows the WRF forecast for $6\textsc{pm}$ for the state of Virginia using the following physics packages (the physics options are given in parentheses):
\begin{itemize}
  \item Micro-physics (Kessler),
   \item Cumulus-physics (Kain),
   \item Short-wave radiation physics (Dudhia),
   \item Long-wave radiation physics (Janjic).
\end{itemize}
Figure \ref{fig:virginia_analysis} shows the NCEP analysis at time $6\textsc{pm}$, which is our proxy for the true state of the atmosphere. The discrepancy between the NCEP analysis and the raw WRF forecast is shown in the Figure \ref{fig:virginia_analysis_err}.
Using the model error prediction we can improve the WRF result by adding the predicted bias  to the WRF forecast. 
The discrepancy between the corrected WRF forecast and the NCEP analysis is shown in the Figure \ref{fig:virginia_improved_err}.  The results show a considerable reduction of model errors as compared to the uncorrected forecast of Figure \ref{fig:virginia_analysis_err}. 
Table \ref{tab:correction_error} shows the minimum and average of original model error vs the improved model errors.

\begin{table}[H]
\begin{center}  
   \begin{tabular}{| l | l | l | l | }
    \hline
   &  $ \textnormal{minimum} (\boldsymbol{\Delta}_{t=6\textsc{pm}}) $  & $\textnormal{average}(\boldsymbol{\Delta}_{t=6 \textsc{pm}})$ 
      \\ \hline
Original forecast   & $6.751 \times 10^{-2}  $ & $  5.025 \times 10^{-1} $ 
     \\ \hline
Improved forecast   & $ 2.134 \times 10^{-4} $  & $ 6.352 \times 10^{-2}  $ 
 \\ \hline
     \end{tabular}
\end{center}
 \caption{The minimum and average of $\boldsymbol{\Delta}_{t=6 \textsc{pm}}$ for the original WRF forecast vs the improved forecast}
  \label{tab:correction_error}
\end{table}
%

\subsection{Experiments for problem one: predicting the norm of precipitation forecast error over the entire domain}
\label{error_model_norm}
We now seek to estimate the two-norm of precipitation model error over the entire continental U.S., which gives a global metric for the accuracy of the WRF forecast, and helps provide insight about the physics configurations that result in more accurate forecasts.
 To this end the following mapping is constructed:

\[
\phi^{\rm error} \left( \Uptheta ,  \| \mathbf{o}_{\tau} \|_2, 
\| \boldsymbol{\Delta}_{\tau} \|_2, 
\| \mathbf{o}_{t=12 \textsc{pm}} \|_2, \bar{\mathbf{o}}_{t=12 \textsc{pm}}
 \right) \approx  \| \boldsymbol{\Delta}_{t=12 \textsc{pm}} \|_2,  \quad   7 \textsc{am} \leq \tau < 12 \textsc{pm}.
\]

To build the training dataset, we run WRF with each of the $252$ different physical configurations. 
The forecast window is 6 hours and the WRF model forecast final simulation time is at $12 \textsc{pm}$. The hourly WRF forecast and discrepancy between
the analysis and WRF forecast is used as training features.

The input features are:
\begin{itemize}
\item  different physics schemes  ($\Uptheta$),
\item the norms of the WRF model predictions at previous time windows,  as well as at the current time
($\| \mathbf{o}_{t=12 \textsc{pm}} \|_2,  \| \mathbf{o}_{\tau} \|_2$, $7 \textsc{am} \leq \tau < 12 \textsc{pm}$), and
\item the norms of past observed discrepancies ($\| \boldsymbol{\Delta}_{\tau} \|_2$,  $7 \textsc{am} \leq \tau < 12 \textsc{pm}$).
\end{itemize}

The output variable is the norm of the discrepancy between WRF precipitation prediction and the NCEP precipitation analysis for the current time window ($ \| \boldsymbol{\Delta}_{t=12 \textsc{pm}} \|_2 $).

We use two different learning algorithms, namely, RF with ten trees in the forest and ANN with four hidden layers, the hyperbolic tangent sigmoid activation function in each layer and linear activation function at last layer.
The number of layers and neurons at each layer is tuned empirically. The total number of samples in the training set is $252$ with $15$ of features.
During the training phase the model learns the effect of interaction of different physical configurations on model error and obtains the approximated function $\widehat{\phi}^{\rm error}$.

In the test phase we feed the approximated function the model information from $1 \textsc{pm}$ to the endpoint of the next forecast window $6 \textsc{pm}$ to predict the norm of the model error  $ \| \widehat{\boldsymbol{\Delta}}_{t=6\textsc{pm}} \|_2$.
\[
\widehat{\phi}^{\rm error} \left( \Uptheta ,  \| \mathbf{o}_{\tau} \|_2, 
\| \boldsymbol{\Delta}_{\tau} \|_2, 
\| \mathbf{o}_{t=6 \textsc{pm}} \|_2, \bar{\mathbf{o}}_{t=6 \textsc{pm}}
 \right) \approx  \| \boldsymbol{\Delta}_{t=6 \textsc{pm}} \|_2,  \quad   1 \textsc{pm} \leq \tau < 6 \textsc{pm}.
\]

\paragraph{Validation of the learned error mapping}
Table \ref{tab:RMSE_error_norm} shows the RMSE between the actual and predicted norms of discrepancies for ANN and RF. The RMSE is taken over the 252 runs with different physics combinations. Both learning models are doing well, with  the ANN giving slightly better results than the RF.
\begin{table}[H]
\begin{center}
    \begin{tabular}{ | l | l | l |}
    \hline
   & $RMSE(\| \widehat{\boldsymbol{\Delta}}_{t=6 \textsc{pm}} \|_2, \| \boldsymbol{\Delta}_{t=6 \textsc{pm}} \|_2 )$  
      \\ \hline
ANN   & $  2.6109\times 10^{-3} $ 
     \\ \hline
RF   & $ 2.9188\times 10^{-3}$ 
 \\ \hline
     \end{tabular}
\end{center}
 \caption{Difference between predicted discrepancy norm $\| \widehat{\boldsymbol{\Delta}}_{t=6 \textsc{pm}} \|_2$ and the reference discrepancy norm $\| \boldsymbol{\Delta}_{t=6 \textsc{pm}} \|_2 $. The $RMSE$ is taken over all test cases.}
  \label{tab:RMSE_error_norm}
\end{table}
%

\paragraph{Analysis of the best combination of physical packages}
Based on our prediction of the norm of model error, the best physics combination that leads to lowest norm of precipitation error
over the entire continental U.S. for the given meteorological conditions is: 
\begin{itemize}
\item the BMJ cumulus parameterization, combined with 
\item the WSM5 micro-physics, 
\item Cam long wave, and 
\item Dudhia short wave radiation physics. 
\end{itemize} 
According to the true model errors, the best physics combination leading to the lowest norm of model error is achieved using the BMJ cumulus parameterization, combined with the WSM5 micro-physics, Cam long wave, and Cam short wave radiation physics.

%
\subsection{Experiments for problem two: identify the physical processes that contribute most to the forecast uncertainty}
\label{physics_configs}
The interaction of different physical processes greatly affects precipitation forecast, and we are interested in identifying the major sources of model errors in WRF. To this end we construct the physics mapping \eqref{eqn:physics_mapping} using the norm and the statistical characteristics of the model-data discrepancy (over the entire U.S.) as input features:
\[
  \phi^{\rm physics} \left( \bar{\boldsymbol{\Delta}}_{t=12 \textsc{pm}}, \| \boldsymbol{\Delta}_{t=12 \textsc{pm}}\|_2  \right) \approx \Uptheta.
\]
Statistical characteristics include the mean, minimum, maximum, and variance of the filed across all grid points over the continental U.S. Note that this is slightly different than \eqref{eqn:physics_mapping} where the inputs are the raw values of these discrepancies for each grid point.
The output variable is the combination of physical processes $\Uptheta$ that leads to model errors consistent with the input pattern $\bar{\boldsymbol{\Delta}}_{t=12 \textsc{pm}}$ and $\| \boldsymbol{\Delta}_{t=12 \textsc{pm}}\|_2$.

To build the dataset, the WRF model is simulated for each of the $252$ different physical configurations, and the mismatches between the WRF forecasts and the NCEP analysis at the end of the current forecast window are obtained. Similar to the previous experiment, the initial conditions used in the WRF model is the NCEP analysis for the May 1st 2017 at $6 \textsc{am}$. The forecast window is 6 hours and the WRF model forecast is obtained for time $12 \textsc{pm}$. The discrepancy between the NCEP analysis at $12 \textsc{pm}$ and WRF forecast at $12 \textsc{pm}$ forms the observable discrepancy for the current forecast window $\boldsymbol{\Delta}_{t=12 \textsc{pm}}$. 
For each of the 252 different physical configurations, this process is repeated and statistical characteristics of the WRF forecast model error  $\bar{\boldsymbol{\Delta}}_{t=12 \textsc{pm}}$, and the norm of model error $\| \boldsymbol{\Delta}_{t=12 \textsc{pm}} \|_2$  are used as  feature values of the function $\phi^{\rm physics}$. 

\paragraph{Validation of the learned physics mapping}
From all the collected data points, $80\%$ (202 samples) are used for training the learning model, and the remaining $20\%$ (50 samples) are used for testing purposes. 

The RF has default ten trees in the forest and ANN has four hidden layers and hyperbolic tangent sigmoid activation function in each layer
with linear activation function at last layer.
The number of layers and neurons at each layer is tuned empirically.
The learning model uses the training dataset to learn the approximate mapping $\widehat{\phi}^{\rm physics}$.
This function is applied to the each of the $50$ test samples $\boldsymbol{\Delta}_{t=12 \textsc{pm}}^{\rm test}$ to obtain the predicted physical combinations $\widehat{\Uptheta}_1$. 
In order to evaluate these predictions, we run the WRF model again with the $\widehat{\Uptheta}_1$ physical setting and obtain the new forecast $\widehat{\mathbf{o}}_{t=12 \textsc{pm}}$, and the corresponding observable discrepancy $\widehat{\boldsymbol{\Delta}}_{t=12 \textsc{pm}}^{\rm test}$.  The RMSE between the norm of actual observable discrepancies and  the norm of predicted discrepancies are shown in Table \ref{tab:performance_physics}. The small values of the difference demonstrates the performance of the learning algorithm.
\begin{table}[H]
\begin{center}
    \begin{tabular}{ | l | l | l |}
    \hline 
   & $ RMSE(\| \widehat{\boldsymbol{\Delta}}_{t=12 \textsc{pm}}^{\rm test} \|_2 , \| \boldsymbol{\Delta}_{t=12 \textsc{pm}}^{\rm test} \|_2) $  
      \\ \hline
ANN   & $ 4.1376  \times 10^{-3} $ 
     \\ \hline
RF   & $ 5.8214 \times 10^{-3} $ 
 \\ \hline
     \end{tabular}
\end{center}
 \caption{The RMSE between estimated discrepancy using predicted physical combinations $ \widehat{\boldsymbol{\Delta}}_{t=12 \textsc{pm}}^{\rm test} $ and the reference discrepancy $\boldsymbol{\Delta}_{t=12 \textsc{pm}}^{\rm test}$.}
  \label{tab:performance_physics}
\end{table}

\paragraph{Analysis of variability in physical settings}
We repeat the test phase for each of the $50$ test samples with the scaled values of observable discrepancies $\boldsymbol{\Delta}_{t=12 \textsc{pm}}^{\rm test}/2$ as inputs, and obtain the predicted physical combinations $\widehat{\Uptheta}_2$. 
Large variability in the predicted physical settings $\widehat{\Uptheta}$ indicate that the respective physical packages variability have a strong influence on the WRF forecast error. We count the number of times the predicted physics $\widehat{\Uptheta}_2$ is
different from $\widehat{\Uptheta}_1$ when the input data spans the entire test data set. 

\begin{figure}[h]
  \begin{centering}
    \includegraphics[width=0.65\textwidth, height=0.48\textwidth]{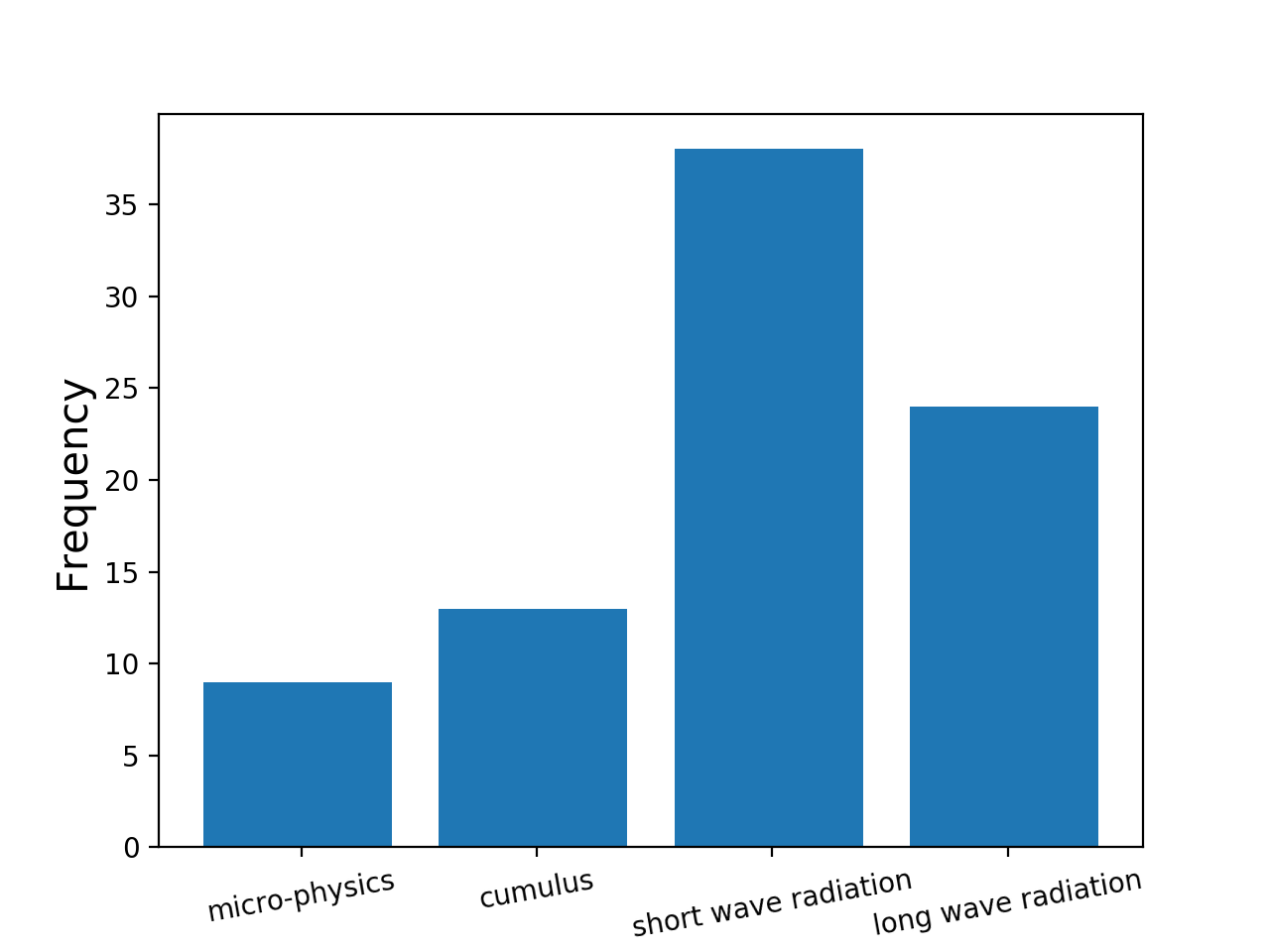}
    \caption{Frequency of change in the physics with respect to change in the input data from $\boldsymbol{\Delta}_{t=12 \textsc{pm}}^{\rm test}$ to $\boldsymbol{\Delta}_{t=12 \textsc{pm}}^{\rm test}/2$. Each data set contains 50 data points, and we report here the number of changes of each package.}
    \label{fig:hist_physics}
  \end{centering}
\end{figure}

The results shown in Figure \ref{fig:hist_physics} indicate that micro-physics and cumulus physics are not too sensitive to the change of input data, whereas short-wave and long-wave radiation physics are quite sensitive to changes in the input data. Therefore our learning model indicates that having an accurate short-wave and long-wave radiation physics package will aid in greatly reducing the uncertainty in precipitation forecasts due to missing/incorrect physics.


%
%
%


\section{Conclusions}
\label{Conc_wrf}
This study proposes a novel use of machine learning techniques to understand, predict, and reduce the uncertainty 
in the WRF model precipitation forecasts due to the interaction of several physical processes included in the model. 

We construct probabilistic approaches to learn the relationships between the configuration of the physical processes used in the simulation and the observed model forecast errors. These relationships are then used to solve two important problems related to model errors, as follows: estimating the systematic model error in a quantity of interest at future times, and identifying the physical processes that contribute most to the forecast uncertainty in a given quantity of interest under specified conditions.

Numerical experiments are carried out with the WRF model using the NCEP analyses as a proxy for the 
real state of the atmosphere. Ensembles of model runs with different parameter configurations are used to generate the training data. Random forests and Artificial neural network models are used to learn the relationships between physical processes and forecast errors. The experiments validate the new approach, and illustrates how it is able to estimate model errors, indicate best model configurations, and pinpoint to those physical packages that influence most the WRF prediction accuracy.

While the numerical experiments are done with WRF, and are focused on forecasting precipitation, the methodology developed herein is general and can be applied to the study of errors in other models, for other quantities of interest, and for learning additional relationships between model physics and model errors.
  \section*{Acknowledgments}
   This work was supported in part by the projects AFOSR DDDAS 15RT1037 and AFOSR Computational Mathematics FA9550-17-1-0205 
  and by the Computational Science Laboratory at Virginia Tech.
The authors would like to thank Dr. R\u{a}zvan \c{S}tef\u{a}nescu for his valuable assistance and suggestions regarding WRF runs and the NCEP dataset.
    %

\section*{References}
\bibliographystyle{plain}
\bibliography{ML_error}

\vfill
\begin{flushright}
\scriptsize \framebox{\parbox{3.2in}{Government License
The submitted manuscript has been created by UChicago Argonne, LLC,
Operator of Argonne National Laboratory (``Argonne").  Argonne, a
U.S. Department of Energy Office of Science laboratory, is operated
under Contract No. DE-AC02-06CH11357.  The U.S. Government retains for
itself, and others acting on its behalf, a paid-up nonexclusive,
irrevocable worldwide license in said article to reproduce, prepare
derivative works, distribute copies to the public, and perform
publicly and display publicly, by or on behalf of the Government. }}
\normalsize 
\end{flushright}
\end{document}